\documentclass[%
 reprint,
 amsmath,amssymb,amsthm,
 aps,
pre,
]{revtex4-2}

\usepackage[version=3]{mhchem}
\usepackage{graphicx,amsmath,amssymb,amsthm}
\usepackage[dvipsnames]{xcolor}
\usepackage{dcolumn}
\usepackage{bm}
\usepackage{bbm}

\def\N{\mathbb{N}}
\def\E{\mathbb{E}}
\def\L{\mathcal{L}}
\newtheorem{theorem}{Theorem}
\newtheorem{defn}{Definition}

\newtheorem{innercustomthm}{Theorem}
\newenvironment{customthm}[1]
  {\renewcommand\theinnercustomthm{#1}\innercustomthm}
  {\endinnercustomthm}

\usepackage{color}
\usepackage[bookmarksnumbered=true]{hyperref} 
\hypersetup{
     colorlinks = true,
     linkcolor = blue,
     anchorcolor = blue,
     citecolor = blue,
     filecolor = blue,
     urlcolor = blue
     }

\definecolor{darkgreen}{rgb}{0.0, 0.5, 0.0}

\newcommand\stdi{\Pi}

\begin{document}

\preprint{APS/123-QED}

\title{Stationary distributions of systems with Discreteness Induced Transitions}

\author{Enrico Bibbona}

 \email{enrico.bibbona@polito.it}
\affiliation{%
Dipartimento di Scienze Matematiche “G.L. Lagrange”\\ Politecnico di Torino, Turin, Italy
}%

\author{Jinsu Kim}
\affiliation{Department of Mathematics
	University of California, Irvine, USA.}
\author{Carsten Wiuf}
\affiliation{Department of Mathematical Sciences\\
	University of Copenhagen, Copenhagen, Denmark
}%

\date{\today}

\begin{abstract}
We provide a theoretical analysis of some autocatalytic reaction networks exhibiting the phenomenon of discreteness induced transitions. The family of networks that we address includes the celebrated Togashi and Kaneko model. We prove positive recurrence, finiteness of all moments, and geometric ergodicity of the  models in the family. For some parameter values, we find the analytic expression for the stationary distribution, and discuss the effect of volume scaling on the stationary behavior of the chain. We  find the exact critical value of the volume for which discreteness induced transitions disappear.
\end{abstract}

\maketitle


\section{
Introduction}

In 2001, Togashi and Kaneko  described a cycle of stochastic  autocatalytic reactions  that displays a highly peculiar dynamics in some regions of the parameter space  \cite{TogashiPRL}. 

It is characterized by switches between patterns where one or more reactants are present in small or vanishing molecule number while other reactants are abundant. The switching is triggered by a single molecule of a previously extinct species that drives the system to a different pattern through  a sequence of quick reactions. The  switches were named Discreteness Induced Transitions (DIT) since  deterministic ODE models are not able to reproduce them \cite{TogashiPRL}. 

The paper raised much interest and  similar effects have been observed in more complicated and realistic models, e.g., large scale networks \cite{largenumsp}, particle  systems with finite interaction radius \cite{vulp1}, reaction-diffusion systems \cite{reacdiff}, models of ant foraging \cite{biancalaniants}, chiral autocatalysis \cite{lentejpc}, tumor growth \cite{cancernib}, spatial models \cite{dipatti}, and viral replication \cite{nibviralrna}.

Several attempts have been made to underpin  the phenomenon theoretically, at least in  simplified toy models, through derivation of analytic expressions, without resorting to simulation or approximation. Examples in this direction are \cite{teoan, anfrem, nibbiancalani,exactnib}, though many questions are still unsolved.

Despite simulation of the Togashi-Kaneko (TK) model indicates  a stationary behavior after a short transient time, positive recurrence (existence of a unique stationary distribution) of the corresponding continuous-time Markov chain (CTMC) has not been proved. For the original 4-dim TK model, no general result from Chemical Reaction Network  theory is applicable. Even if the system is reduced to dimension two,  the problem of finding a stationary distribution remains non-trivial, and the curious switching behavior persists.
In dimensions four and two, the switching behavior  causes  the \emph{seemingly} stationary distribution  emerging from simulation to be multimodal for certain  parameters values. When the rates are scaled in the volume $V$ of the container and $V$ is considered  large, the multimodality disappears and a distribution with a concentrated peak emerges. In this case, the scaled stochastic model converges to the classical deterministic model (fluid limit). 

In this paper, we prove that a family of autocatalytic networks, including  the TK model, is positive recurrent in arbitrary dimension (Theorem \ref{thm:non-explosive}). For some parameter values, an explicit expression for the stationary distribution is derived. In 2-dim (cf. Theorem \ref{th:mixed}) the parameter region for which  the stationary distribution is known covers the 2-dim TK model. In higher dimension (cf. Theorem \ref{th:mixed2})  the parameter region for which  the stationary distribution is known, does not include the general TK model. However, it includes a large family of TK-like  models exhibiting DITs..

The analytic form of the stationary distribution provides a clear theoretical demonstration of the effect of volume scaling on the stationary behavior of the system. It also allows us to find the exact critical value of the volume from which the DIT stops to appear. This value we  also find for the TK model.

\section{Background material}

\subsection{
The original 4-dim TK model}

Let $\N=\{0,1,2,\ldots\}$ denote the integers including zero.
For any two integers $i$ and $n$, let $( i)_n$ be the remainder after integer division of $i$ by $n$ (elsewhere denoted by $i \text{ mod }n$).
The network proposed by Togashi and Kaneko  \cite{TogashiPRL} consists of the following cycle of autocatalytic reactions
\begin{equation} A_i+A_{({i+1})_4}\ce{->[\kappa]}2 A_{({i+1})_4},\qquad i=1,\ldots,4, \label{ciclo1} \end{equation}
together with inflow and outflow reactions
\begin{equation*}
	A_i \ce{<=>[\delta][\lambda]}\emptyset,\qquad i=1,\ldots,4.
\end{equation*}
The state of the system is a tuple of  four non-negative integers $\mathbf{a}=(a_1,a_2,a_3,a_4)'$. 
Denote by $\mathbf{e}_j$  the $j$-th unit vector,  $j=1,\ldots,4$. The transitions rates generated by the autocatalytic reactions are
\[q_{\mathbf{a},\mathbf{a}-\mathbf{e}_i+ \mathbf{e}_{({i+1})_4}}=\kappa a_i a_{({i+1})_4},\]
while those corresponding to inflow and outflow reactions are
\[q_{\mathbf{a},\mathbf{a}+\mathbf{e}_i}=\lambda\quad \text{and}\quad q_{\mathbf{a},\mathbf{a}-\mathbf{e}_i}=\delta a_i. \]
The qualitative behavior of the system depends on the parameter values. The classical volume scaling (cf.~\cite[Chapter 11]{kurtzbook} or \cite{bibbonareview}) is adopted in \cite{TogashiPRL},  where the initial molecule counts of the species are proportional to the  scaling parameter $V$. It implies the  rate constants are given by 
\[\kappa=\frac{\kappa'}{V}\qquad \delta=\delta' \qquad \lambda=\lambda'V\]  

One parameter can always be set to one by linear scaling of time. In \cite{TogashiPRL},  $\kappa'=1$, and further $\lambda'=\delta'=D$ for simplification.
According to \cite{kurtz1ode} or \cite[Chapter 11]{kurtzbook}, when $V\rightarrow\infty$,  the density process, which is the CTMC rescaled by dividing the molecule numbers by $V$, converges to the solution of a system of ordinary differential equations with stable equilibrium $(1,1,1,1)$. Indeed when  $VD\gg1$, the reaction rates  are large and the trajectories of the density process only display small fluctuations around the deterministic equilibrium. 

For  $VD\ll1$, a completely different behavior appears,  triggered by the slow rate of inflow and DIT appears. If the system is initialized at a state where all species counts are large, one of the species at random (say, species 3) is quickly driven to extinction by the fast autocatalytic dynamics. At this point, several molecules of species 2 are produced and not consumed and they catalyze the consumption of all molecules of species 1. 
We end up with a configuration where the species 1 and 3 are both consumed, the count of species 2 is very high, and that of species 4 is quite low. We call this pattern 2H4L.
In this configuration only slow inflows and outflows are active, and one needs to wait until a molecule of species 3 or  1 flows in before the autocatalytic dynamics starts again leading to  another pattern with two non-contiguous species extinct.
The dynamics of the system then proceeds by switching between such patterns in a way that a 2H4L configuration is much more often followed by a 2L4H pattern and only rarely switches to either a 1H3L or 1L3H configuration (cf. FIG 1 in  \cite{TogashiPRL}).  

\subsection{
Lumpability}

In the next Section we exploit the notion of lumpability to find the stationary distribution in some cases. We summarize here the meaning of this property.

Let $\{S_I\}_{I\in \mathcal I}$ be a partition of a denumerable state space $S$ of a CTMC $X(t)$, $t\ge 0$, with rates $q_{ij}$, $i,j\in S$.
Let moreover $\iota$ be the function that maps  $x\in S$ to the index of the element of the partition to which $x$ belong, (i.e $\iota(x)=K $ if and only if $x\in S_K$).
The process $X(t)$, $t\ge 0$,  is (strongly) lumpable if the \emph{lumped} process $\overline X(t)=\iota(X(t))$, $t\ge 0$,  is a CTMC on $\mathcal I$ for any choice of initial distribution. Sufficient conditions (cf. \cite{lumpabilityCTMC}) that guarantees lumpability of a regular, irreducible, positive recurrent CTMC $X(t)$, $t\ge 0$,  on the partition $\{S_I\}_{I\in \mathcal I}$ are that  every subset $S_I$ is finite, and that for any $I,J\in \mathcal{I}$, and any $i,i'\in S_I$, 
\[\sum_{j \in S_J} q_{ij}=\sum_{j \in S_J} q_{i'j}=\overline{q}_{IJ} \]

The rates of the lumped chain $\overline X(t)$, $t\ge 0$,  are  $\overline{q}_{IJ}$, $I,J\in\mathcal{I}$, and for any $s< t$, the lumped variable $\overline X(t)$ is independent of $X_s$ given $\overline X_s$

\section{The 2-dim TK model}

If the number of species in the TK model is reduced to two, the reaction network becomes 
\begin{equation} \label{ciclo2d}
\begin{aligned} 2A_1\ce{<-[\kappa_1]} A_1&+A_2\ce{->[\kappa_2]} 2 A_{2} \\
A_1 \ce{<=>[\delta_1][\lambda_1]}&\hspace{1 mm}\emptyset\hspace{1 mm}\ce{<=>[\lambda_2][\delta_2]}A_2,
\end{aligned}
\end{equation}
where we allow $\kappa_1$ and $\kappa_2$ to be different.
The state of the network is denoted by $\mathbf{a}=(a_1, a_2)'\in \N^2$, the molecule counts of each species.
The transitions rates due to the autocatalytic reactions are
\begin{equation}q_{\mathbf{a},\mathbf{a}-\mathbf{e}_1+ \mathbf{e}_2 }=\kappa_2 a_1 a_2, \qquad q_{\mathbf{a},\mathbf{a}-\mathbf{e}_2+ \mathbf{e}_1}=\kappa_1 a_1 a_2,\label{acrates}\end{equation}
while those corresponding to inflow and outflow reactions are
\begin{equation}q_{\mathbf{a},\mathbf{a}+\mathbf{e}_i}=\lambda_i, \qquad q_{\mathbf{a},\mathbf{a}-\mathbf{e}_i}=\delta_i a_i,  \qquad  i=1,2.
\label{otherrates}\end{equation}

\begin{figure}
	\includegraphics[width=7.5 cm]{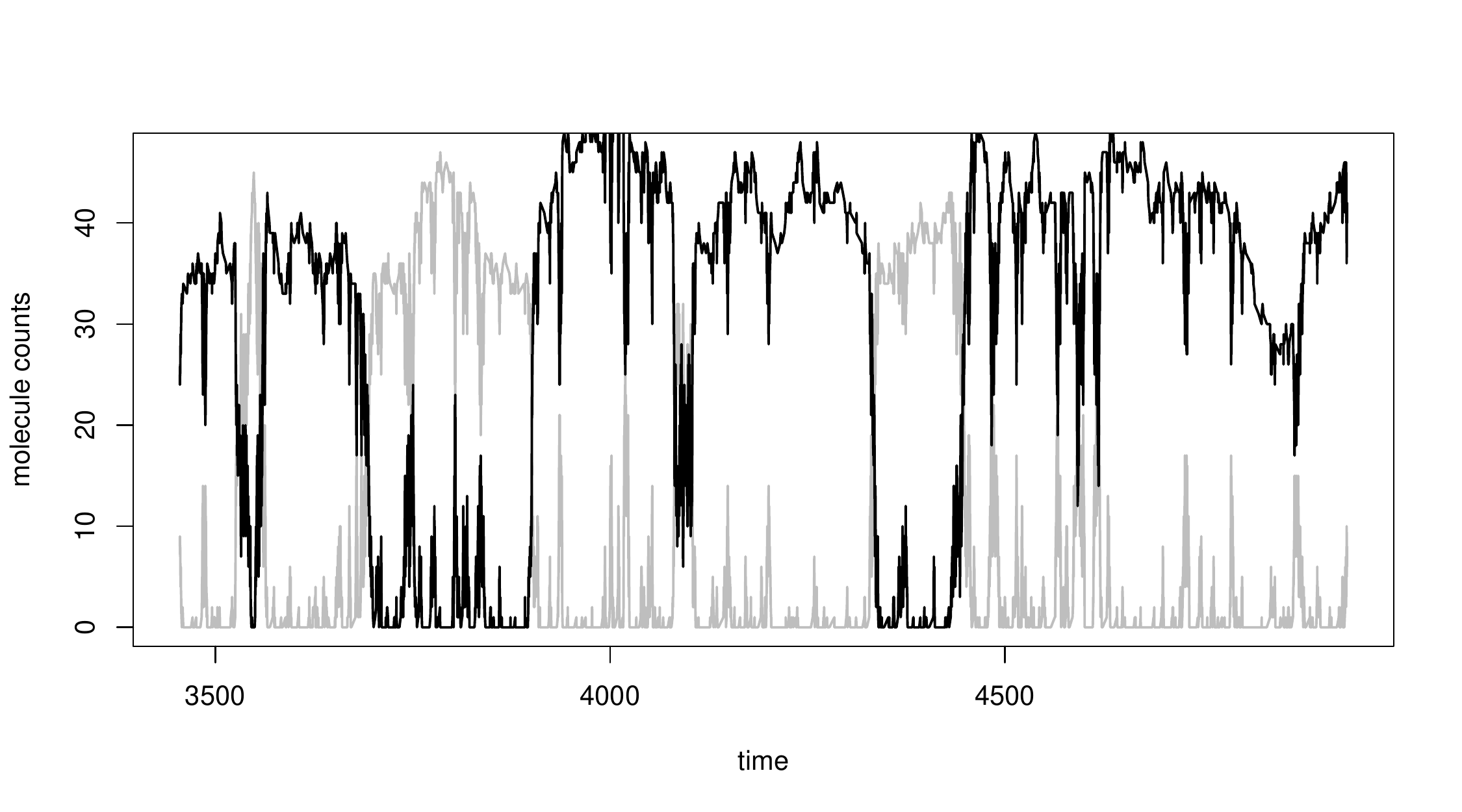}
	\caption{\label{fig1}Molecule counts of the two species of network \eqref{ciclo2d} along time. Patterns where the grey species is mostly absent alternates with patterns where the  mostly absent species is the black one. Parameters are $\lambda_i=0.2$, $\delta_i=0.01$, $\kappa_i=0.05$ for  $i=1,2$.}
\end{figure}

The dynamics is  simplified, but not too much. When the inflows occur at a much slower rate than the autocatalytic reactions, the system switches between two patterns in a similar way to the original $4$-dimensional TK system, where one or the other compound is mostly absent. 
A plot of the two simulated trajectories in this parameter range  is shown in FIG. \ref{fig1}.

\subsection{Positive recurrence and stationary distribution}

A proof of positive recurrence in a  more general setting is given in Section \ref{subsec:positive recurrence for d}. In this section, we show the sketch of how to derive the stationary distribution of \eqref{ciclo2d} by using its lumpability. The subsets 
\[E_n=\{\mathbf{a} \in \N^2 \colon a_1+a_2=n\},\]
 form a partition $\{E_n\}_{n\in\N}$ of the state space. The CTMC model $X(t)$, $t\ge 0$,  of \eqref{ciclo2d} under stochastic mass-action kinetics is lumpable with respect to this partition if $\delta=\delta_1=\delta_2$.
 With this choice, the rate at which the total molecule count $n$ is increased by one, is equal to the sum of the rates of the inflows
\begin{equation*}q_{n,n+1}=\sum_{i=1}^2 q_{\mathbf{a},\mathbf{a}+\mathbf{e}_i}=\lambda_1+\lambda_2, \label{rate+}\end{equation*}
independently of  $\mathbf{a}$.
The rate at which  $n$ is decreased by one, is the sum of the rates of the outflows
 \begin{equation*}q_{n,n-1}=\sum_{i=1}^2 q_{\mathbf{a},\mathbf{a}-\mathbf{e}_i}=\delta(a_1+a_2)=\delta n\label{rate-}\end{equation*}
and therefore it does not depend on $\mathbf{a}$ as long as $\mathbf{a}\in E_n$.

The lumped process $\overline X(t)$, $t\ge 0$,  is described by the following reaction network where a single species $B$ aggregates all  molecules of $A_1$ and $A_2$
\begin{equation}B\ce{<=>[\lambda_1+\lambda_2][\delta]}\emptyset. \label{cb}\end{equation}

Network \eqref{cb} is weakly reversible and has  \emph{deficiency}  zero \cite{andersonKurtzBook}. By   \cite[Theorem 3.6 and 3.7]{andersonKurtzBook}, it admits a unique stationary distribution with Poisson law 
 \begin{equation}\label{piconn}
 \nu(n)=\frac{\mu^n}{n!}\exp\left(-\mu\right), \quad \mu=\frac{\lambda_1+\lambda_2}{\delta}.
 \end{equation}
where $n$ is the state of the lumped process (i.e., $\overline X(t)=n$ if and only if $X(t)\in E_n$). 
We now aim at factorizing the stationary distribution $\stdi(\mathbf{a})$ of the process $X(t)$, $t\ge 0$, of  \eqref{ciclo2d} by conditioning on the stationary probability $\nu(n)$ of the lumped process being in state $n=a_1 +a_2$. We write
   \begin{equation}
   \stdi(\mathbf{a})=\pi(a_1|n)\nu(n).
   \label{factorized}\end{equation}
   
A careful rewriting of the master equation for the stationary distribution $\stdi(\mathbf{a})$ shows that $\stdi(\mathbf{a})$ is stationary if and only if $\pi(\mathbf{a}|n)$ fulfils
\begin{equation}R_n= L_{n-1}+L_{n}+L_{n+1}\label{eq:pi}\end{equation}
with
\begin{align}
R_n=&(\lambda_1+\lambda_2+ n\delta +(\kappa_1+\kappa_2)a(n-a))\pi (a|n)\nonumber \\
L_{n-1}=&  \frac{n\delta\lambda_1}{\lambda_1+\lambda_2} \pi(a-1| n-1) +\frac{n\delta \lambda_2}{\lambda_1+\lambda_2}\pi(a| n-1) \nonumber \\
L_{n}=& \kappa_1(a+1)(n-a-1) \pi(a+1|n)+\nonumber \\
&+\kappa_2 (a-1)(n-a+1) \pi(a-1|n) \nonumber \\
L_{n+1}=&\frac{\lambda_1+\lambda_2}{n+1}(a+1)\pi(a+1|n+1)\nonumber\\ &+\frac{\lambda_1+\lambda_2}{n+1}(n- a+1)\pi(a|n+1), \nonumber
\end{align}
for $n\ge 0$ and $a=0,\ldots,n$.

\begin{figure}
	\includegraphics[width=0.5 \textwidth]{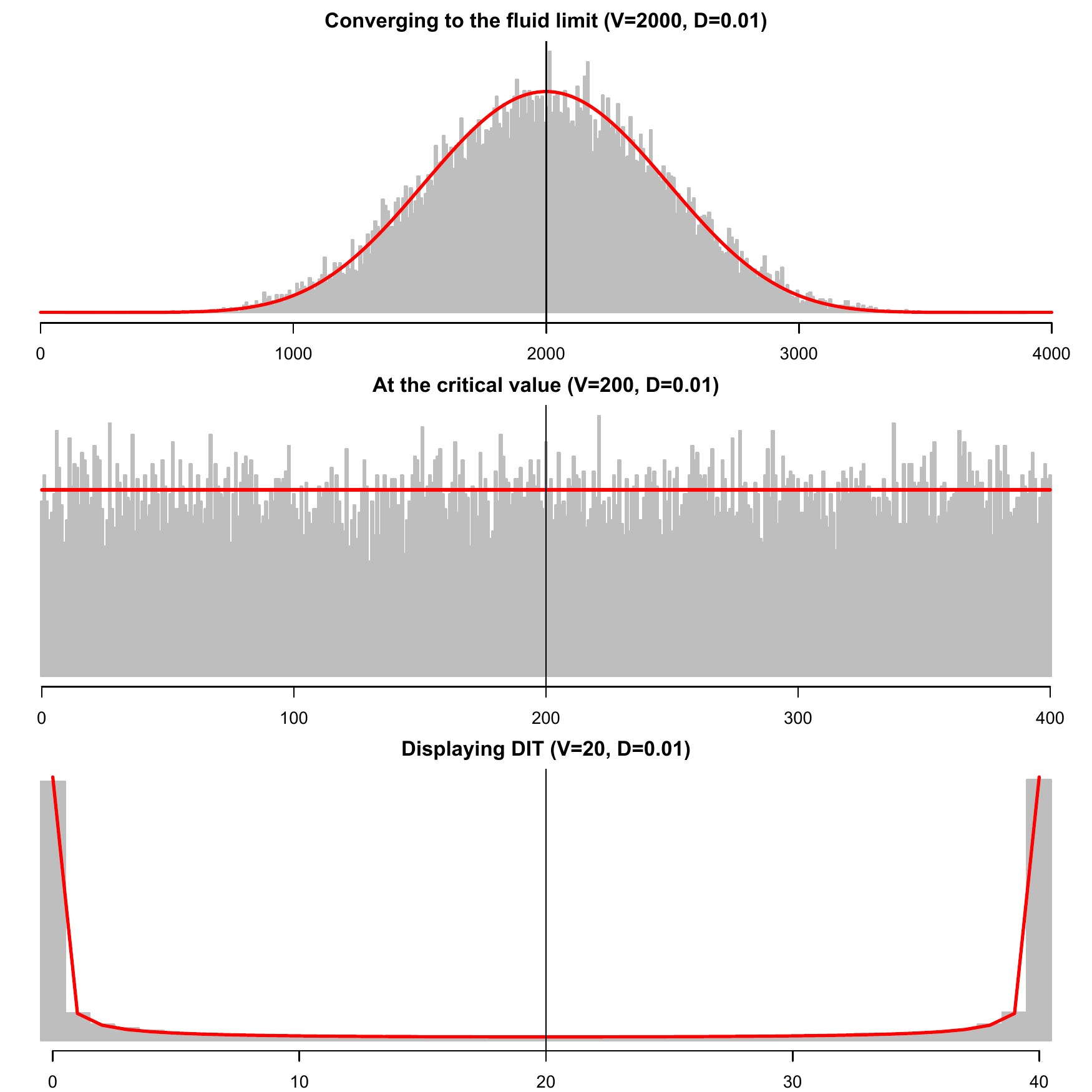}
	\caption{\label{fig3} 
	The effect of scaling and DIT. The conditional stationary distribution $\pi(a_1|n)$, cf. \eqref{factorized} from simulation (grey histograms) and from  \eqref{bebi} in Theorem \ref{th:mixed} (red lines).  Parameters are chosen according to \eqref{scaledconstants} with $\kappa'_i=1$, $\delta'_i=\lambda'_i=D=0.01$ for $i=1,2$.   
	The volume parameter $V$ differs in the three panels to illustrate the effect of  scaling, and $n$ is chosen as the mean of $\nu(n)$, which is 4000, 400, and 40, respectively; implying that the mean of the scaled process $X/V$ is $(1,1)$  in all three cases. 
	Simulation set-up: $2.75\cdot 10^6$, $10^6$, and $10^5$ simulations (from upper to lower panel) were conducted with fixed time $T=250, 50, 50$, respectively (the stationary regime already applies).  Only values of $\mathbf{a}(T)$ with $a_1(T)+a_2(T)=n$ were kept, and the histogram of $a_1(T)$ was plotted.}
\end{figure}

Unfortunately, there is not a simple way to find a closed form expression of $\pi(\cdot |n)$ satisfying  equation \eqref{eq:pi}.
  However, simulation of the process for different  rate constants, corresponding to different regimes of the volume  $V$ (cf. Section \ref{VS} for more details), indicates that the conditional stationary distribution may  be unimodal, flat or concentrated at the boundaries (cf. FIG \ref{fig3}).  Statistical practice suggests the beta-binomial as a natural candidate for a discrete distribution on the integers $\{0,\cdots, n\}$ that may display these behaviors. The next theorem confirms this, and FIG \ref{fig3} provides a graphical comparison between simulations and theoretical values in different parameter settings.

\begin{theorem}\label{th:mixed} 
	Network \eqref{ciclo2d}, assuming that $\kappa=\kappa_1=\kappa_2>0$ and $\delta=\delta_1=\delta_2>0$, has  a unique stationary distribution $\stdi(\mathbf{a})$ that factorizes as \eqref{factorized}, where $\nu(n)$ is  given by \eqref{piconn}, and $\pi(\cdot|n)$ is given by the beta-binomial distribution
	\begin{equation}\pi(i|n)=\binom{n}{i}\frac{B(i+\alpha,n-i+\beta)}{B(\alpha,\beta)},\quad i=0,\ldots,n,\label{bebi}\end{equation}
	where 
\begin{equation}\label{alphabeta}\alpha=\frac{\delta\lambda_1}{\kappa(\lambda_1+\lambda_2)},\quad \beta=\frac{\delta\lambda_2}{\kappa(\lambda_1+\lambda_2)},\end{equation}
	and
	$$B(x,y)=\frac{\Gamma(x)\Gamma(y)}{\Gamma(x+y)},\quad x,y>0.$$
\end{theorem}
\begin{proof}
	The proof is by direct verification, substituting expression \eqref{bebi} into equation \eqref{eq:pi}. Calculations are displayed  in Appendix \hyperref[app:stationary]{B} in a more general context.
	\end{proof}

\subsection{Volume scaling}\label{VS}

Molecule counts and mass-action rates  can be scaled with the volume $V$ in such a way that the scaled stochastic system $X(t)/V$ converges for large $V$ to the solution of the deterministic system on any finite time horizon, cf. \cite{kurtzbook}[Ch. 11, Theorem 2.1].

This is achieved for \eqref{acrates} and \eqref{otherrates}, under the hypothesis of Theorem \ref{th:mixed}, by setting the constants to
\begin{equation}
\kappa_i=\frac{\kappa'}{V}\quad\delta_i=\delta' ,\quad \lambda_i=\lambda'_i V,
\label{scaledconstants}\end{equation}
for  $i=1,2$.
When $V$ is not sufficiently large  the stochastic model differs significantly from the deterministic limit \cite{TogashiPRL} and starts to display the switching behavior (DIT)  illustrated in FIG \ref{fig1}.

\begin{figure}
	\includegraphics[width=0.5 \textwidth]{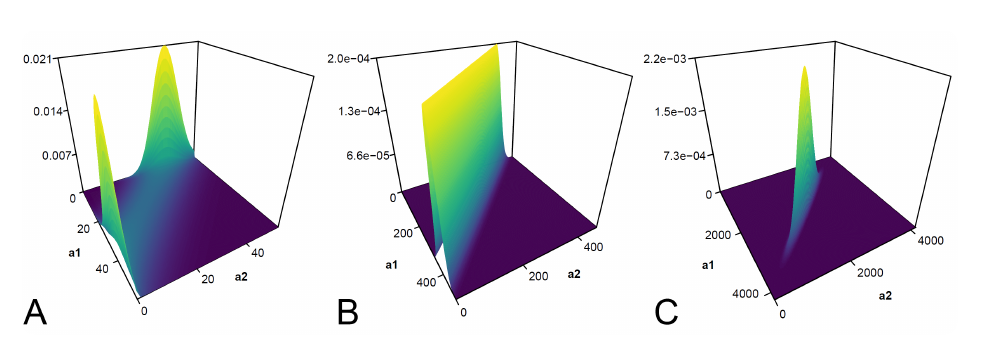}
	\caption{\label{fig2}Smoothed representation of the stationary density \eqref{factorized}. In all panels $D=0.01$. The volume parameter $V$ is different in the three panels to illustrate the effect of the scaling. In panel A, $V=20$ (same range as in FIG. 1). The density is bimodal and concentrated at the boundaries, since DIT are present. In panel B, $V=200$ and the conditional density \eqref{bebi} is uniform. In panel C, $V=2000$ and the density is concentrated around the deterministic equilibrium.}
\end{figure}

In \cite{TogashiPRL}, the authors set 
\begin{equation}\kappa'_i=1, \qquad \delta'_i=\lambda'_i=D,\qquad i=1,2.\label{scaledconstants1}\end{equation}
With this choice of the rate constants, by Theorem \ref{th:mixed}, we know the explicit form of the stationary density, and we can investigate the behavior of the system at every $V$ without resorting to simulations. 
The stationary conditional density $\pi(\cdot|n)$ is beta-binomial with parameters $\alpha=\beta=DV/2$.
The beta-binomial density is unimodal when $\alpha$ and $\beta$ are both larger than one (that is, when  $DV$ is larger than two)  with the mass concentrated at  the equilibrium of the corresponding deterministic model. When both $\alpha$ and $\beta$ are smaller than one
(that is, when $DV$ is small) the density becomes bimodal with most of the mass at the boundaries. The intermediate case is when $\alpha=\beta=1$ and the conditional distribution reduces to the discrete uniform distribution on $\{0,\ldots,n\}$. In other words, at the critical value $DV=2$, the conditional density flattens to $\pi(a|n)=\frac{1}{n+1}$ for every $a$. A pictorial representation of the density \eqref{factorized}, at different values of $V$ with $D$ fixed to the value $0.01$ is given in FIG. \ref{fig2}. The effect of the scaling is apparent. For graphical convenience, the discrete density has been smoothed to a continuous one.

To make this effect quantitatively apparent, in the general setting where \eqref{scaledconstants} holds but not necessarily  \eqref{scaledconstants1}, we prove that for $V\rightarrow 0$ the stationary distribution concentrates at the boundaries by showing that the conditional probability  $\pi(0|n)+\pi(n|n)$ tends to one, for any $n$. 
Indeed, inserting \eqref{scaledconstants} into \eqref{alphabeta}, we get $\alpha=\alpha' V$ and $\beta=\beta' V$ with \[\quad\alpha'=\frac{\delta\lambda'_1}{\kappa(\lambda'_1+\lambda'_2)} \quad\text{and} \quad \beta'=\frac{\delta\lambda'_2}{\kappa(\lambda'_1+\lambda'_2)}. \]
The sum of the two conditional probabilities  
reduces to 
\begin{align}
 {}&\pi(0|n)+\pi(n|n)=\label{corners}\\
 &=\left[\frac{\Gamma(n+\alpha' V)]}{\Gamma(\alpha' V)}+\frac{\Gamma(n+\beta' V)}{\Gamma(\beta' V)}\right]\frac{\Gamma[(\alpha'+\beta')V}{\Gamma[n+(\alpha'+\beta')V]}.\notag
 \end{align}
Whatever  $n$ is, since $\Gamma(z)\sim 1/z$ for $z\rightarrow 0$, it is easily seen that the sum tends to one as $V\rightarrow 0$.

For large $V$, we show that the stationary distribution $\stdi_V$ of the scaled process $X(t)/V$ concentrates around the deterministic equilibrium $(1,1)$. The mean $\mu_V$ and variance $\Sigma_V$ of $\stdi_V$ might easily be computed (by conditioning on $n$) from the first and second moments of the Poisson distribution  and  the beta-binomial distribution.

The explicit calculation is here only reported for two components, but can be found  for the others as well,
\begin{align*}
\left(\mu_V\right)_1&= \frac{\mu}{V}\frac {\alpha}{\alpha+\beta}\\
\left(\Sigma_V\right)_{11}&=\frac{1}{V^2}\left[\frac {\alpha\beta}{(\alpha+\beta)^{2}}\frac {(\alpha+\beta) \mu+ \mu^2+\mu}{\alpha+\beta+1} + \mu \frac {\alpha}{\alpha+\beta}\right], \end{align*}
where $\mu$ is given in \eqref{piconn} and $\alpha$ and $\beta$ in \eqref{alphabeta}.
Scaling the parameters as  in  \eqref{scaledconstants}, it is easily observed that $\left(\mu_V\right)_1 \rightarrow 1$  and $\left(\Sigma_V\right)_{11} \rightarrow 0$ for $V\rightarrow \infty$.
With a little more effort,  the same result  extends to the other components, that is, we have
\[\mu_V \to \left(\frac{\lambda'_1}{\delta'}, \frac{\lambda'_2}{\delta'}\right) \quad\text{ and }\quad \Sigma_V \to \begin{pmatrix}0&0\\ 0&0\end{pmatrix}.\]

In general, the \emph{agreement} between the stochastic and the deterministic model for large volumes only holds on a finite time horizon only. Negative examples where the two modelling paradigms differ asymptotically are discussed in \cite{cappellettijmb,agazzimattingly}. Our result shows that for large $V$, under the assumptions of Theorem \ref{th:mixed}, the stochastic and the deterministic models of \eqref{ciclo2d} are in agreement  asymptotically.

\section{
Higher dimensional models}\label{sec4}

In higher dimension there exist different models whose 2-dim reduction corresponds to network \eqref{ciclo2d}. One of them is the 4-dim TK model \eqref{ciclo1}, but also the network
\begin{equation}
2A_i \ce{<-[\kappa_{ji}]}  A_i+A_{j}\ce{->[\kappa_{ij}]} 2 A_j
\qquad A_i \ce{<=>[\delta_i][\lambda_i]}\emptyset,
\label{4dv3}
\end{equation}
$i,j=1,\ldots,d$, $i\neq j$, can be seen as a $d$-dimensional version of model \eqref{ciclo2d}. Network \eqref{4dv3} includes \eqref{ciclo1} as a special case for $\kappa_{ij}$ equal to $\kappa$ when $j=(i+1)_d$ and zero otherwise.
Reaction rates are the obvious generalizations of \eqref{acrates} and \eqref{otherrates}.

\subsection{Positive recurrence}\label{subsec:positive recurrence for d}

In this section, we state the positive recurrence of the Markov process underlying the general $d$-dimensional model \eqref{4dv3}. To do so, we show that $V(x)=e^{\Vert x \Vert_1}$, where  $\Vert x \Vert_1=\sum_{i=1}^d|x_i|$,  is a Lyapunov function. Non-explosivity and positivity, then, follow by the Foster-Lyapunov criterion \cite{MT-LyaFosterIII}. Additionally, as a by-product, all moments of the stationary distribution are positive and convergence to the stationary distribution is exponentially fast.
The detailed proof  can be found in Appendix~\hyperref[app:Positive recur]{A}.
 
\begin{theorem}\label{thm:non-explosive}
For any non-negative values of the parameters $\kappa_{ij}$, $\kappa_{ji}$, and for positive $\lambda_i$ and $\delta_{i}$, the CTMC  associated  to the  system \eqref{4dv3}  is positive recurrent on $\N^d$ (for any $d$). Consequently, it has  a unique stationary distribution supported on $\N^d$. Moreover, all moments are finite and the convergence to the stationary distribution is exponentially fast.
\end{theorem}

\subsection{Stationary distribution}

\subsubsection{The model  }

By the same argument as we used in dimension 2, under the assumption of equal outflow rates ($\delta_i=\delta$ for all $i= 1\cdots d$), the process $X(t)$, $t\ge 0$,  that counts the molecules of each species is lumpable on the partition $\{E_n\}$, where $E_n=\{\mathbf{a} \in \N^2\colon  \sum_i a_i=n\}$.

The lumped process $\overline X (t) =\sum_i X_i(t)$ represents the total molecule count. It follows a birth and death process (as in  \eqref{cb}) with Poisson stationary distribution  with intensity
\begin{equation}\mu=\frac{\sum_{i=1}^d \lambda_i}{\delta}.\label{newlambda1}\end{equation}
Similarly to the 2-dim case,  the stationary distribution $\stdi(\mathbf{a})$ factorizes as
\begin{equation}\stdi(\mathbf{a})=\pi(\mathbf{a}|n)\nu(n).\label{factorized12}\end{equation}

\begin{theorem}\label{th:mixed2} 
Assume $\kappa_{ij}=\kappa>0$, $i,j=1,\ldots,d$, $i\not=j$,  $\delta=\delta_1=\cdots =\delta_d >0$, and $\lambda_i >0$ for all $i$. Then, model \eqref{4dv3} has a unique stationary distribution $\stdi(\mathbf{a})$ expressed as in \eqref{factorized12}, where $\nu(n)$ is given as in \eqref{piconn} and \eqref{newlambda1}, and $\pi(\cdot|n)$ is given by the Dirichlet-multinomial distribution
	\begin{equation}
	\pi(\mathbf{a}|n)=\binom {n}{\mathbf{a}}\frac{\Gamma(\sum_{i=1}^d \alpha_i)}{\Gamma(n+\sum_{i=1}^d \alpha_i)}\prod_{i=1}^d\frac{\Gamma(a_i+\alpha_i)}{\Gamma(\alpha_i)}\label{DM11}\end{equation}
	where $\mathbf{a}$  is any $d$-dimensional integer vector  with  $\Vert \mathbf{a} \Vert_1=n$,  and 
	$$\alpha_i=\frac{\delta\lambda_i}{\kappa \sum_{i=1}^d \lambda_i}.$$
\end{theorem}

\begin{proof}
	The proof is by direct verification, substituting expression \eqref{DM11} into equation \eqref{factorized12} using \eqref{piconn}.
	Calculations are displayed  in Appendix \hyperref[app:stationary]{B}.
\end{proof}

\subsection{Volume scaling and other properties}

The scaled process $X(t)/V$ in dimension $d$ has similar properties to that of the scaled process in dimension 2. 
  In the case where the stationary distribution is known (Theorem \ref{th:mixed2}), we might proceed similarly to what was done in dimension 2 and calculate the mean vector and covariance matrix of the molecule counts, now using moment properties of the Poisson and the Dirichlet-multinomial distributions. 
  Parameters are scaled according to
  \begin{equation}
  \kappa=\frac{\kappa'}{V}\quad\delta=\delta' ,\quad \lambda_i=\lambda'_iV, \label{parametriscalati}
 \end{equation}
  for  $i=1,\cdots, d$.
  
  As the volume $V$ increases towards infinity, it can be shown that the mean vector converges to $\left(\tfrac{\lambda_1}{\delta},\ldots,\tfrac{\lambda_d}{\delta}\right)$, the equilibrium point of the deterministic process, and the covariance matrix decreases towards the zero matrix. Thus, under the hypothesis of Theorem \ref{th:mixed2}, the deterministic and the stochastic models of \eqref{4dv3} are in agreement asymptotically for large volume size in the long run (at stationarity) as well as over the finite time horizon.

At the other extreme, for $V\to 0$,  the conditional probability of a corner configuration tends to one.
Indeed, such probability, that generalizes \eqref{corners}, is equal to
 \[\sum_{i=1}^d \pi(\mathbf{e}_i|n)= \frac{ \Gamma(V\sum_j \alpha'_j)}{\Gamma(n+V\sum_j \alpha'_j)}\sum_{i=1}^d\frac{\Gamma(n+V\alpha'_i)}{\Gamma(V\alpha'_i)},\]
 where \[\alpha'_i=\frac{\delta\lambda'_i}{\kappa \sum_j\lambda'_j}.\]
The convergence to one can be easily shown with the same methodology used in dimension 2.

\begin{figure}
	\includegraphics[width=0.5 \textwidth]{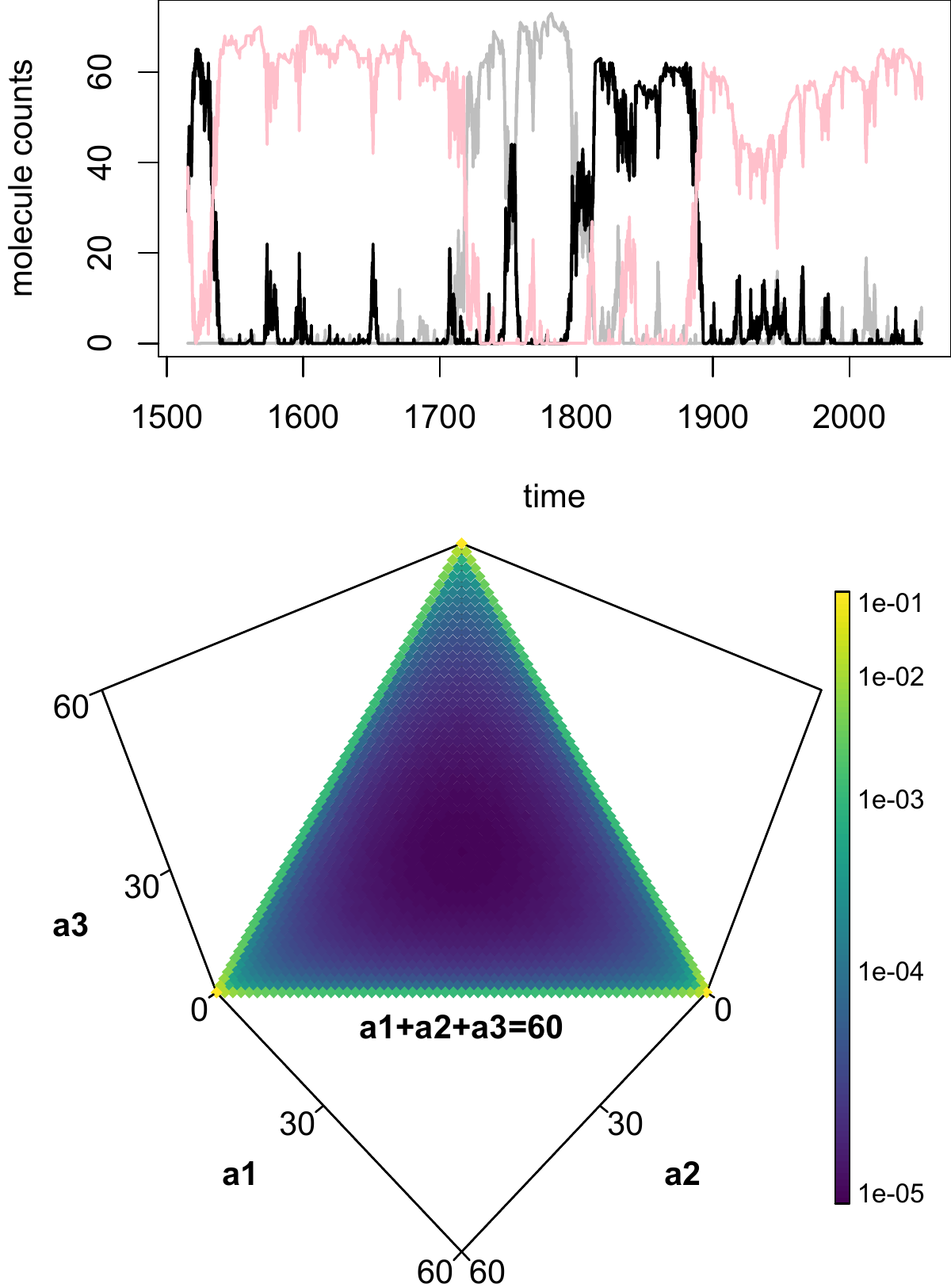}
	\caption{\label{fig4} Simulated trajectories and the conditional stationary distribution $\pi(\mathbf{a}|n)$, from model \eqref{4dv3} in dimension three.  Parameters are chosen according to \eqref{volll} with $D=0.01$ and $V=20$.  The value of $n$ is fixed to 60 in the lower plot. The presence of DITs is apparent both from the trajectories and from the conditional distribution that is concentrated at the corners of the simplex.}
\end{figure}

The peaks at the vertexes reflect the presence of DIT that causes the switch between dynamical patterns where only one of the species is present in large quantity at a time, while all the others are almost extinct.
A graphical illustration of the presence of DITs in a three dimensional version of model \eqref{4dv3} is given in FIG. \ref{fig4}. In dimension three it is no longer possible to plot the stationary distribution $\stdi(\mathbf{a})$. However, choosing $D=0.01$ and $V=20$, we can plot a set of simulated trajectories, and the values of the conditional stationary distribution $\pi(\mathbf{a}|n)$.

If parameter are further chosen as
\begin{equation}\kappa'_i=1, \qquad \delta'_i=\lambda'_i=D,\qquad i=1,\cdots, d \label{volll} ,\end{equation}
in analogy of what was done in \cite{TogashiPRL}, the distribution becomes symmetric in the labels of the species
and the $\alpha_i=DV/d$, $i =1,\ldots,d$, are all equal. Moreover if $V=\frac{d}{D}$, the conditional distribution is flat, providing a transition  point from the multimodal case to the unimodal case. If the $\alpha_i$s are not equal (i.e. the $\lambda'_i$ are not), the transition will not proceed through a flat conditional distribution.

Other relevant properties of the Dirichlet-multinomial distribution, like aggregation, marginals, conditional distributions are discussed in \cite{mosimann62,hoadley69}.

\subsection{Back to the d-dim TK model}

Model \eqref{ciclo1} motivated our interest in autocatalytic networks. Theorem \ref{thm:non-explosive} guarantees that it is positive recurrence, but an explicit form of the stationary distribution cannot be derived by Theorem \ref{th:mixed2}. Indeed, it is a special case of \eqref{4dv3}, where some of the $\kappa_{ij}$ are  set to zero (those for which $j \neq (i+1)_d$) and all others are set to the same value $\kappa$.  However, it is still possible to find the explicit expression in a very special case.
\begin{theorem}\label{th:tk} 
	Assume that $\kappa=\kappa_1=\cdots= \kappa_d\ge 0$ and $\delta=\delta_1=\cdots =\delta_d=\frac{d}{d-1}\kappa$ and $\lambda=\lambda_1=\cdots =\lambda_d>0$. Then, model \eqref{ciclo1} has a unique stationary distribution $\stdi(\mathbf{a})$ whose expression is \eqref{factorized12} with $\nu(n)$ given by \eqref{piconn} and \eqref{newlambda1} and with $\pi(\cdot|n)$ given by a uniform distribution
	 \begin{equation}\pi(\mathbf{a}|n)=\frac{n!(d-1)!}{(n+d-1)!}\label{unif1}
	 \end{equation} on the simplex $\{\mathbf{a}\in \{0,\ldots n\}^d \colon\Vert \mathbf{a} \Vert_1=n\}$.
\end{theorem}
	The proof is by direct verification, substituting expression \eqref{unif1} into equation \eqref{factorized12} with  $\kappa_{ij}$  set to zero for all $j \neq (i+1)_d$ and  to the same value $\kappa$ otherwise. Calculations are displayed in Appendix \hyperref[app:TK]{C}. If the rate constants are scaled in the volume as in equation \eqref{parametriscalati} and further set to \eqref{volll}, the critical value of the volume that makes the distribution flat is $V=\frac{d}{(d-1)D}$, in agreement with the result for $d=2$. In 4-dim, in \cite{TogashiPRL}, it was noticed from simulation that the order of the magnitude of this critical value should be around $V\sim 1/D$. However,  determining the exact value was pursued.  Our result allows us to ensure that the exact value is $V=4/3D$.

\begin{acknowledgments}
This paper started with a group discussion at the Erwin Schr\"{o}dinger Institute (ESI) workshop on “Advances in Chemical Reaction Network Theory”, Vienna (15-19/10, 2018) and was further developed during the Summer school and Workshop on Chemical reaction networks organised at Politecnico di Torino (24/6-3/7, 2019). We acknowledge the organizers and the ESI for the invitation and the financial support. EB gratefully acknowledges funding from the Italian Ministery of Education, University and Research, MIUR, grant Dipartimenti di Eccellenza 2018-2022 (E11G18000350001).  JK thankfully acknowledges  travel support from NSF grant DMS1616233 to German Enciso. Computational resources for the preliminary simulation were provided by HPC@POLITO (http://www.hpc.polito.it).

\end{acknowledgments}

\section*{Appendix}\label{sec:nd}
In this Appendix, we report the detailed proofs of the three main theorems  (Theorem \ref{thm:non-explosive}, \ref{th:mixed2}, \ref{th:tk}) that are stated in Section \ref{sec4}. Theorem \ref{th:mixed} is not proved separately since it is a special case of Theorem \ref{th:mixed2}.
The most general high dimensional model that we consider is that in \eqref{4dv3}, that we repeat here
\begin{equation}
2A_i \ce{<-[\kappa_{ji}]}  A_i+A_{j}\ce{->[\kappa_{ij}]} 2 A_j
\qquad A_i \ce{<=>[\delta_i][\lambda_i]}\emptyset.\notag
\end{equation}

\subsection*{Appendix A} 

As we stated in the main text, the associated CTMC for the general dimensional model \eqref{4dv3} is positive recurrent and admits a unique stationary distribution. We prove this in the following theorem. We further show that the CTMC is \textit{exponentially ergodic} meaning that the associated distribution $P^t$ at time $t$ converges to the unique stationary distribution exponentially fast. The proof relies on the Foster-Lyapunov criterion \cite{MT-LyaFosterIII}. We begin with a formal statement and necessary concepts for the Foster-Lyapunov criterion.

\begin{defn}
For a CTMC $X(t)$, $t\ge 0$, defined on a countable state space $\chi$, the infinitesimal generator $\mathcal{L}$ is the operator  
\begin{align*}
\mathcal L V(x)=\sum_{\eta} \lambda_\eta (x)(V(x+\eta)-V(x)),
\end{align*}
 where $\eta$ is a transition of $X(t)$, $\lambda_\eta$ is the transition rate associated with $\eta$, and $V$ is any real function defined on the state space.
\end{defn}

For a CTMC $X(t)$, $t\ge 0$, we define a truncated process $X_M$ such that $X_M(t)=X(t)$ if $|X(t)| < M$ and $X_M(t)=x_M$ otherwise for some fixed state $x_M$ with $|x_M|\ge M$. We denote by $\mathcal{L}_M$ the infinitesimal generator of $X_M$. 

We further call $V(x)$ a norm-like function if $V(x)$ is a positive function such that $|V(x)|\to \infty$, as $|x|\to \infty$.

The following theorem is Theorem 6.1 in \cite{MT-LyaFosterIII}, in the case of a countable state space. It  is one version of the Foster-Lyapunov criterion for exponential ergodicity.

\begin{customthm}{A}[Foster-Lyapunov criterion \cite{MT-LyaFosterIII}]\label{thm:fl}
Let $X(t)$, $t\ge 0$, be a CTMC defined on a countable state space $\chi$.  Then $X(t)$, $t\ge 0$, is non-explosive and positive recurrent if there exist a norm-like function $V$ on $\chi$, positive constants $C$ and $D$ such that for any $M>0$

\begin{equation*}
\mathcal{L}_M V(x) \le -CV(x)+D \quad \text{for all $x \in \chi$}.
\end{equation*}
Furthermore, $X(t)$, $t\ge 0$, admits a unique stationary distribution $\pi$ on each irreducible component, and there exist $B>0$ and $\beta \in (0,1)$ such that
\begin{align*}
\sup_{A} |P^t(x,A)-\pi(A)| \le B V(x) \beta ^t \text{ for all $x\in \chi$}.
\end{align*}
\end{customthm}

To show  positive recurrence and exponential ergodicity of the CTMC associated with the general model \eqref{4dv3}, it is therefore sufficient to prove that there exist a norm-like function $V$ and positive constants $C, D$ such that 
\begin{equation}\label{eq:the goal}
\mathcal{L}V(x) \le -CV(x)+D \quad \text{for all $x$}.
\end{equation}
In the proof of the following theorem, we prove \eqref{eq:the goal} for an exponential function $V$. Using this specific function, we also show that all moments of the unique stationary distribution of $X(t)$, $t\ge 0$, are finite.

\setcounter{theorem}{1}

\begin{proof}[Proof of Theorem \ref{thm:non-explosive}]
Let $X(t)$, $t\ge 0$, be the CTMC associated with the system \eqref{4dv3}.
	 Let $V(x)=e^{\Vert x \Vert_1}$, where  $\Vert x \Vert_1=\sum_{i=1}^d|x_i|$. Then we show that \eqref{eq:the goal} holds for some positive constants $C$ and $D$.
	
	Let $e_i\in\N^d$ be the vector with $i$-th component  $1$ and zero otherwise. We have
	\begin{align*}
	\mathcal{L}V(x)=&\sum_{i,j} \kappa_{ij}x_ix_j(V(x+e_i-e_j)-V(x))\\
	&+ \sum_{i,j} \kappa_{ij}x_ix_j(V(x+e_i-e_{j})-V(x))\\
	&+\sum_{i=1}^d \delta_{i}x_i(V(x-e_i)-V(x))\\
	&+\sum_{i=1}^d \lambda_{i}(V(x+e_i)-V(x))\\
	\end{align*}
	\begin{align*}
	\mathcal{L}V(x)=&\sum_{i=1}^d \delta_{i}x_i(V(x-e_i)-V(x))\\
	&+\sum_{i=1}^d \lambda_{i}x_i(V(x+e_i)-V(x)).
	\end{align*}
	
	Let $K_n=\{x \in \N^d \colon x_i \ge n \text{ for each  i}\}$. Then note that for $x\in K_n$,
	\begin{align*}
	\mathcal{L}V(x)
	& = V(x)\left (\sum_{i=1}^d\delta_ix_{i}(e^{-1}-1)+\sum_{i=1}^d\lambda_i(e-1)\right)\notag\\
	& \le \left(\min_i\delta_i(e^{-1}-1)dn+\sum_{i=1}^d\lambda_i(e-1)\right)V(x) 
	\end{align*}
	Hence, by choosing sufficiently large $N$ such that $$C=-\left(\min_i\delta_i(e^{-1}-1)dN+\sum_{i=1}^d\lambda_i(e-1)\right)>0,$$ we conclude that \eqref{eq:the goal} holds with $D=2C\max_{x\in K^c_N}V(x)$. This implies that $X(t)$, $t\ge 0$, is non-explosive, positive recurrent and exponential ergodicity by Theorem \ref{thm:fl}. This implies existence of a unique stationary distribution $\pi$.

	 To show that $\pi$ has finite $m$th moment for any $m\in\N^d$,  we use \eqref{eq:dynkin2} below combined with the ergodic theorem \cite{NorrisMC97}. Let $\tau_M=\inf \{ t>0:|X(t)| \ge M \}$. Then by using Dynkin's formula \cite{dynkin1965markov, oksendal2013stochastic} and \eqref{eq:the goal}, we have 
	\begin{equation*} 
     \begin{split}	
	&\E_{x}(V(X(t)))=V(x)+\E_{x}\!\left (\int_0^{t\wedge \tau_M} \mathcal{L} V(X(s)) ds \right) \quad \ \ \notag\\
	&\le V(x)-C\E_{x}\left( \int_0^{{t\wedge \tau_M}} V(X(s))ds \right ) + Dt \notag \\
	\end{split}
    \end{equation*}
    \vspace{-0.5cm}		
    \begin{equation}
	\hspace{0.3cm} =V(x)-C\E_{x}\!\left( \int_0^t V(X(s)) \mathbbm{1}_{\{|X(s)|<M\}} ds \right ) + Dt \label{eq:dynkin2},
\end{equation}			
	where $\E_{x}$ denotes the expectation of $X(t)$ with $X(0)=x$, and $t\wedge \tau_M =\min\{t,\tau_M\}$.
	 By rearranging terms in \eqref{eq:dynkin2} and dividing by $t, C$, it follows that
	\begin{equation}\label{eq:dynkin3}
	\E_{x}\!\left(\frac{1}{t}\int_0^{t} V(X(s))\mathbbm{1}_{\{|X(s)|<M\}} ds \right ) \le \frac{V(x)}{Ct}+\frac{D}{C}.
	\end{equation}
	 
	Then by the dominant convergence theorem, taking $\lim$ for $t\to \infty$ on both sides in \eqref{eq:dynkin3} gives that
	\begin{align*}
	&\lim_{t \to \infty}\E_{x}\!\left( \frac{1}{t}\int_0^t  V(X(s))\mathbbm{1}_{\{|X(s)|<M\}}ds \right )\\
	&=\sum_{x\in \N^d}V(x)\mathbbm{1}_{\{|x|<M\}} \pi(x)\le \frac{D}{C}.
	\end{align*}
	Then the monotone convergence theorem applies for $M\to \infty$ to conclude that $\sum_{x\in \N^d}V(x)\pi(x)\le \frac{D}{C}$.
	Since $V(x)=e^{\Vert x \Vert_1}$, any moment of $\pi$ is finite.
\end{proof}

\subsection*{Appendix B. \ Stationary distribution}\label{app:stationary}

\begin{proof}[Proof of Theorem \ref{th:mixed2}]
	Under the assumption of equal outflow rates, the process $X(t)$ that counts the molecules of each species is lumpable on the partition $\{E_n\}_{n\in\N}$, where $E_n=\{\mathbf{a} \in \N^2\colon  \sum_{i=1}^d a_i=n\}$.
	
	The lumped process $\overline X (t) =\sum_{i=1}^d X_i(t)$ 
	has Poisson stationary distribution $\nu(n)$  with intensity \eqref{newlambda1}.
	As stated earlier, the stationary distribution $\stdi(\mathbf{a})$ factorizes as $\stdi(\mathbf{a})=\pi(\mathbf{a}|n)\nu(n)$.
	Under the given assumptions on the parameters, $\pi(\mathbf{a}|n)$ solves the  equation, similar to \eqref{eq:pi},

	\begin{equation}R_n= L_{n-1}+L_{n}+L_{n+1},\label{eq:pi2}\end{equation}
	where 
	\begin{align}
	R_n=&\pi (\mathbf{a}|n) \left[\sum_{i=1}^d \lambda_i + \delta n +\sum_{i=1}^d \sum_{j\neq i}\kappa a_i a_j \right]\nonumber \\
	L_{n-1}=&  \frac{\delta n}{\sum_{i=1}^d \lambda_i}\sum_{i=1}^n \pi(\mathbf{a}-\mathbf{e}_i| n-1) \lambda_i \nonumber \\
	L_{n}=& \sum_{i,j=1}^n \pi(\mathbf{a} -\mathbf{e}_i + \mathbf{e}_j| n) \kappa (a_i-1) (a_j+1)\nonumber \\
	L_{n+1}=&\frac{\sum_{i=1}^n \lambda_i}{(n+1) }\sum_{i=1}^n \pi(\mathbf{a} +\mathbf{e}_i |n+1) (a_i+1),\notag
	\end{align}
	
	The proof now proceeds by showing that the ansatz $\pi(\cdot|n)$ specified by equation \eqref{DM11} solves equation \eqref{eq:pi2}.
	First we note that if the ansatz is true, then the following recurrence relations hold
	\begin{equation} 
	\begin{split}
	&\pi(\mathbf{a}|n)=\frac{1}{n+1}\sum_{i=1}^d (a_i +1)\, \pi(\mathbf{a}+\mathbf{e}_i|n+1),\\
	&\pi(\mathbf{a}-\mathbf{e}_i|n-1)= \frac{a_i (n-1+\sum_{i=1}^d \alpha_i)}{n(a_i-1+\alpha_i)}\pi(\mathbf{a}|n),\\
	&\pi(\mathbf{a}-\mathbf{e}_i+\mathbf{e}_j|n)=  \frac{a_i(a_j +\alpha_j)}{(a_j+1)(a_i-1+\alpha_i)}\pi(\mathbf{a}|n).
     \end{split}\label{recurrence}
	\end{equation}
	Applying \eqref{recurrence} and dividing by $\pi (\mathbf{a}|n)$ in \eqref{eq:pi2} we get
	\begin{align}\label{sospeso}	
	&\delta n +\sum_{i=1}^d\sum_{j\neq i}\kappa a_i a_j 	=\frac{\delta  (n-1+\sum_{i=1}^d \alpha_i)}{\sum_{i=1}^d \lambda_i}\sum_{i=1}^d  \frac{\lambda_i a_i }{a_i-1+\alpha_i} \notag \\
	&+\sum_{i=1}^d \sum_{j\neq i} \frac{\kappa a_i (a_i-1) (a_j +\alpha_j)}{a_i-1+\alpha_i} 		
	\end{align}
	By fixing $a_i=n$, the following condition is necessary
	\begin{equation}
	\begin{split}
	&\delta (n-1+\alpha_i) \\
	&=\delta  \Big(n-1+\sum_{i=1}^d \alpha_i\Big)\frac{\lambda_i }{\sum_{i=1}^d \lambda_i}+  \kappa (n-1) \sum_{j\neq i}\alpha_j. \label{augualen}
    \end{split}	
	\end{equation}
If we further set $n=1$ we get 
	\begin{equation}
	\frac{\lambda_i}{{\sum_{i=1}^d \lambda_i}} =\frac{\alpha_i}{\sum_{i=1}^d \alpha_i}.\label{primo}
	\end{equation} 
	Moreover if we take equation \eqref{augualen} and sum over all $i=1,\cdots, d$, we get
	\[(d-1)\delta (n-1)+\delta \sum_{i=1}^d \alpha_i = \delta\sum_{i=1}^d \alpha_i +  \kappa (n-1)(d-1) \sum_{i=1}^d \alpha_i,\]
	which further implies
	\begin{equation}
	\sum_{i=1}^d \alpha_i  =\frac{\delta}{\kappa}.\label{sumalphas}
	\end{equation}
Together with equation \eqref{primo}, this implies
	\begin{equation} 
	\alpha_i= \frac{\delta\lambda_i}{{\kappa\sum_{i=1}^d \lambda_i}}.\label{alphas}
	\end{equation}
	
	Taking again equation \eqref{sospeso}, we can further recast it into the following form
	\begin{align*}\delta n=& -\sum_{i=1}^d\sum_{j\neq i}\kappa a_i a_j \notag+\frac{\delta  (n-1+\sum_{i=1}^d \alpha_i)}{\sum_{i=1}^d \lambda_i}\sum_{i=1}^d  \frac{\lambda_i a_i }{a_i-1+\alpha_i}\notag\\
	&+\sum_{i=1}^d \sum_{j\neq i} \frac{\kappa a_i (a_i-1+\alpha_i) (a_j +\alpha_j)}{a_i-1+\alpha_i} \notag\\
	&-\sum_{i=1}^d \sum_{j\neq i} \frac{\kappa a_i \alpha_i(a_j +\alpha_j)}{a_i-1+\alpha_i} \notag\\
	=&\frac{\delta  (n-1+\sum_{i=1}^d \alpha_i)}{\sum_{i=1}^d \lambda_i}\sum_{i=1}^d  \frac{\lambda_i a_i }{a_i-1+\alpha_i}+\notag\\
	&+\sum_{i=1}^d \sum_{j\neq i} \kappa a_i  \alpha_j -\sum_{i=1}^d \sum_{j\neq i} \frac{\kappa a_i \alpha_i(a_j +\alpha_j)}{a_i-1+\alpha_i} \notag\\
	=&\frac{\delta  (n-1+\sum_{i=1}^d \alpha_i)}{\sum_{i=1}^d \lambda_i}\sum_{i=1}^d  \frac{\lambda_i a_i }{a_i-1+\alpha_i}+ \notag\\
	&+\sum_{i=1}^d \sum_{j=1}^d \kappa a_i  \alpha_j - \sum_{i=1}^d  \kappa a_i  \alpha_i \notag\\
	&-\sum_{i=1}^d \sum_{j=1}^d \frac{\kappa a_i \alpha_i(a_j +\alpha_j)}{a_i-1+\alpha_i}+
	\sum_{i=1}^d \sum_{i=1}^d \frac{\kappa a_i \alpha_i(a_i +\alpha_i)}{a_i-1+\alpha_i}, \notag\\
	\end{align*}
\begin{align*}
	\delta n  =&\frac{\delta  (n-1+\sum_{i=1}^d \alpha_i)}{\sum_{i=1}^d \lambda_i}\sum_{i=1}^d  \frac{\lambda_i a_i }{a_i-1+\alpha_i}+\sum_{i=1}^d \sum_{j=1}^d \kappa a_i  \alpha_j \\
	&- \sum_{i=1}^d  \kappa a_i  \alpha_i -\Big(n +\sum_{i=1}^d\alpha_i\Big)\sum_{i=1}^d \frac{\kappa a_i \alpha_i}{a_i-1+\alpha_i}\\
	&+	\sum_{i=1}^d \sum_{i=1}^d \frac{\kappa a_i \alpha_i(a_i +\alpha_i)}{a_i-1+\alpha_i} \notag\\
	=&\frac{\delta  (n-1+\sum_{i=1}^d \alpha_i)}{\sum_{i=1}^d \lambda_i}\sum_{i=1}^d  \frac{\lambda_i a_i }{a_i-1+\alpha_i} + \sum_{i=1}^d \sum_{j=1}^d \kappa a_i  \alpha_j\\  
	&- \sum_{i=1}^d  \kappa a_i  \alpha_i -\Big(n +\sum_{i=1}^d\alpha_i\Big)\sum_{i=1}^d \frac{\kappa a_i \alpha_i}{a_i-1+\alpha_i}\\
	&+	\sum_{i=1}^d \sum_{i=1}^d \frac{\kappa a_i \alpha_i(a_i +\alpha_i-1)}{a_i-1+\alpha_i}\notag\\
	&-\sum_{i=1}^d \sum_{i=1}^d \frac{\kappa a_i \alpha_i}{a_i-1+\alpha_i} \notag
	\end{align*}
	and, finally,

	\begin{align}\delta n  =&\frac{\delta  (n-1+\sum_{i=1}^d \alpha_i)}{\sum_{i=1}^d \lambda_i}\sum_{i=1}^d  \frac{\lambda_i a_i }{a_i-1+\alpha_i}	 \label{ultima}\\
	+&\sum_{i=1}^d \sum_{j=1}^d \kappa a_i  \alpha_j -\Big(n +\sum_{i=1}^d\alpha_i-1\Big)\sum_{i=1}^d \frac{\kappa a_i \alpha_i}{a_i-1+\alpha_i}\notag
     \end{align}

	Now, using \eqref{sumalphas} and \eqref{alphas}, we have
	\[\sum_{i=1}^d  \kappa a_i  \sum_{j=1}^d\alpha_j= n \delta\]
	and
	\[\frac{\delta}{\sum_{i=1}^d \lambda_i}\sum_{i=1}^d  \frac{\lambda_i a_i }{a_i-1+\alpha_i}=\sum_{i=1}^d \frac{\kappa a_i \alpha_i}{a_i-1+\alpha_i},\]
	making equation \eqref{ultima} identically satisfied.
	\end{proof}

\subsection*{Appendix C} 

\begin{proof}[Proof of Theorem \ref{th:tk}.]
We first rewrite equation \eqref{eq:pi2}, where the  $\kappa_{ij}$ are set to zero for all $j \neq (i+1)_d$ and  to a constant value $\kappa$ otherwise. We get the condition
\begin{equation}R_n= L_{n-1}+L_{n}+L_{n+1},\label{eq:pi3}\end{equation}
where 
\begin{align}
R_n=&\pi (\mathbf{a}|n) \left[\sum_{i=1}^d \lambda_i + \delta n +\sum_{i=1}^d \kappa a_i a_{(i+1)_d}\right]\nonumber \\
L_{n-1}=&  \frac{\delta n}{\sum_{i=1}^d \lambda_i}\sum_{i=1}^n \pi(\mathbf{a}-\mathbf{e}_i| n-1) \lambda_i \nonumber \end{align}
\begin{align}
L_{n}=& \sum_{i=1}^n \pi(\mathbf{a} -\mathbf{e}_i + \mathbf{e}_{(i+1)_d}| n) \kappa (a_i-1) (a_{(i+1)_d}+1)\nonumber \\
L_{n+1}=&\frac{\sum_{i=1}^n \lambda_i}{(n+1) }\sum_{i=1}^n \pi(\mathbf{a} +\mathbf{e}_i |n+1) (a_i+1).\notag
\end{align}
We now notice that if the uniform ansatz is true, the following recurrence relations also hold
\[\begin{aligned} \pi(\mathbf{a}+\mathbf{e}_i|n+1)=&\frac{n+1}{n+d}\: \pi(\mathbf{a}|n) \\ \pi(\mathbf{a}-\mathbf{e}_i|n-1)=&\frac{n+d-1}{n}\: \pi(\mathbf{a}|n)
\end{aligned}\]
Plugging the ansatz \eqref{unif1} and these recurrence relations into \eqref{eq:pi3}, we get that equation \eqref{eq:pi3} holds if and only if
\[\delta n +\sum_{i=1}^d\kappa a_i a_{(i+1)_d} =\delta  (n +d-1)+\sum_{i=1}^d  \kappa (a_i-1) (a_{(i+1)_d} +1).\]
It simplifies to $0=\delta (d-1)-\kappa d.$
Such a condition is identically satisfied under the hypothesis of the theorem which guarantees
\[\kappa=\frac{d-1}{d}\delta.\]
\end{proof}

\bibliography{bibl.bib}
\end{document}